\RequirePackage{fix-cm}
\documentclass[smallextended, envcountsame]{svjour3}
\smartqed
\def\makeheadbox{\relax}
\journalname{}

\usepackage[unicode = false,
    pdftoolbar = true,
    colorlinks = true,
    linkcolor = blue,
    citecolor = blue,
    filecolor = black,
    urlcolor = blue,
    breaklinks = true]{hyperref}
\usepackage[ruled]{algorithm2e}
\usepackage{float}
\usepackage[utf8]{inputenc}
\usepackage{listings}           
\usepackage{amsmath,amssymb}            
\usepackage{xcolor}             
\usepackage{graphicx}           
\usepackage[ruled]{algorithm2e}
\usepackage{verbatim}
\usepackage{geometry}
\usepackage{enumerate}
\usepackage[thinlines, thiklines]{easybmat}
\usepackage{hhline}
\usepackage{nicematrix}
\usepackage{mathtools}
\usepackage[noblocks]{authblk}

\usepackage{graphicx}
\usepackage{pgfplots}
\usepackage{tikz}
\usetikzlibrary{patterns}
\usetikzlibrary{arrows}
\usetikzlibrary{shapes}
\pgfplotsset{compat=newest}
\usepgfplotslibrary{fillbetween}

\DeclareMathOperator{\spann}{span}
\DeclareMathOperator{\pspan}{pspan}

\newcommand{\R}{\mathbb{R}}

\newcommand{\cm}[1]{\operatorname{cm}\left(#1\right)}

\newcommand{\V}{\mathrm{cV}}

\newcommand{\bbm}{\begin{bmatrix}}
\newcommand{\ebm}{\end{bmatrix}}
\newcommand{\G}{\mathbf{G}}

\newcommand{\zero}{\mathbf{0}}

\newcommand{\Proj}{\operatorname{Proj}}
\newcommand{\ProjD}{\operatorname{P}_D}
\newcommand{\ProjV}{\operatorname{P}_V}
\newcommand{\dom}{\operatorname{dom}}

\newcommand{\cl}{\operatorname{cl}}

\newcommand{\N}{\operatorname{\mathbb{N}}}

\newcommand{\A}{\operatorname{\mathcal{A}}}

\newcommand{\one}{\mathbf{1}}

\newcommand{\AL}{\overline{\A}}
\DeclareMathOperator{\proj}{Proj}

\newcommand{\p}{\mathring{p}_{B}}
\newcommand{\us}{u_{B}}
\DeclareMathOperator{\CM}{cm}
\DeclareMathOperator{\CV}{cV}

\newcommand{\CMD}{\CM_D(D)}
\newcommand{\CVD}{\V_D(D)}
\DeclareMathOperator{\argminop}{argmin}
\newcommand{\argmin}[1]{\underset{#1}{\argminop}}

\definecolor{gris9}{rgb}{0.9,0.9,0.9}
\definecolor{gris8}{rgb}{0.8,0.8,0.8}
\definecolor{gris7}{rgb}{0.7,0.7,0.7}
\definecolor{gris6}{rgb}{0.6,0.6,0.6}
\definecolor{gris5}{rgb}{0.5,0.5,0.5}
\definecolor{gris4}{rgb}{0.4,0.4,0.4}
\definecolor{gris3}{rgb}{0.3,0.3,0.3}
\definecolor{gris2}{rgb}{0.2,0.2,0.2}
\definecolor{gris1}{rgb}{0.1,0.1,0.1}

\begin{document}

\title{\centering The cosine measure relative to a subspace}
\author{\centering  Charles Audet \and Warren Hare \and Gabriel Jarry-Bolduc}
\institute{Charles Audet \at
            {GERAD}
		and D\'epartement de math\'ematiques et g\'enie industriel,
		Polytechnique Montr\'eal,
		Montr\'eal, Qu\'ebec, Canada. 
            Audet's research is partially funded by the Natural Sciences and Engineering Research Council (NSERC) of Canada, Discovery Grant \#2020-04448.
            ORCID 0000-0002-3043-5393 \newline
		\email{charles.audet@polymtl.ca}
\and
Warren Hare \at
              Department of Mathematics, University of British Columbia, Kelowna, British Columbia, Canada. \\ Hare's research is partially funded by the Natural Sciences and Engineering Research Council (NSERC) of Canada, Discover Grant \#2018-03865. ORCID 0000-0002-4240-3903\\
              \email{warren.hare@ubc.ca}           
           \and
          Gabriel Jarry-Bolduc\at
           Department of Mathematics and Statistics, Saint Francis Xavier University, Antigonish, Nova Scotia, Canada. \\ ORCID 0000-0002-1827-8508\\
              \email{gabjarry@alumni.ubc.ca} 
              }
\date{\today}

\titlerunning{The cosine measure relative to a subspace}

\maketitle
\begin{abstract}
The \emph{cosine measure} was introduced in 2003 to quantify the richness of finite positive spanning sets of directions in the context of derivative-free directional methods.
A positive spanning set is a set of vectors whose nonnegative linear combinations span the whole space.
The present work extends the definition of cosine measure.
In particular, the paper studies cosine measures relative to a subspace, 
    and proposes a deterministic algorithm to compute it.
The paper also studies the situation in which the set of vectors is infinite.
The extended definition of the cosine measure might be useful for subspace decomposition methods.
\end{abstract}

\section{Introduction}\label{sec:intro}

Given a set of vectors $D$ in $\R^n$, the concepts of {\em spanning}, {\em linear independence}, and {\em basis}, are considered foundational to linear algebra. Closely related, but less well studied, are the ideas of {\em positive spanning}, {\em positive linear independence}, and {\em positive basis}~\cite{Davis1954}.  All of these notions can be defined considering a set's properties relative to $\R^n$ or relative to a linear subspace \cite{regis2016}.

In addition to the intrinsic mathematical interest in these concepts, {positive bases} have been shown to be fundamental to the convergence analysis in {\em direct-search methods} in {\em derivative-free optimization}.
The first occurrence dates from 1996~\cite{RMLewis_VTorczon_1996}
    and more recent presentations are found in the textbooks~\cite{audet2017derivative} and~\cite{conn2009introduction}.
A key tool in this analysis is the {\em cosine measure}, introduced in~\cite{kolda2003}.  Algorithms to compute the cosine measure  are provided in \cite{hare2020deterministic,regis2021}. Results on  the maximal value of the cosine measure for positive bases of $\R^n$ are developed in \cite{hare2023nicely,naevdal2018}. The notion of cosine measure is also  briefly investigated for \emph{positive $k$-spanning set} of $\R^n$ in \cite{hare2024using}.
Until now, this cosine measure was only deeply studied for positive spanning sets of $\R^n$.  In this paper, we explore the {\em cosine measure relative to a subspace}, which effectively provides a measure of how well a subset of the positive spanning set considered explores that subspace.  The results herein will be of high value in understanding the {\em subspace decomposition methods} which are recently gaining traction in \cite{AuDeLe08,CartisRobert2023,Gratton2017direct,kimiaei2023subspace,kozak2021stochastic,xie2023new}.

The main goals of this paper are to introduce the notion of cosine measure relative to a subspace, to demonstrate its value in understanding positive spanning properties of a set, and to provide a deterministic algorithm to compute it.  As a secondary goal, we provide several novel results examining the case where the set of vectors is infinite.  

The remainder of this paper is organized as follows. 
The notation is presented in Section~\ref{sec:prel} and  background results that are necessary to understand this paper are provided.  This includes {\em positive spanning}, {\em positive linear independence}, {\em positive basis},  {\em cosine measure}, {\em cosine vector set}, and {\em active set}.  These final three definitions are extended to include ``relative to a subspace''.  In Section \ref{sec:properties}, properties of a positive spanning set of a linear subspace are investigated. It is shown that many known results regarding $\R^n$ are easily adapted to working in a subspace or working with an infinite set.  In Section \ref{sec:CMrelativetoL}, the  notion of cosine measure relative to a subspace is explored.  Several results are provided that link the positive span, cosine measure, and cosine measure relative to a subspace.  The case where the cosine measure relative to a subspace is equal to $0$ is examined.  The section concludes with new results showing how the cosine measure relative to a subspace can be used to provide a general error bound on the true gradient of a smooth function. In Section \ref{sec:algo}, a deterministic algorithm to compute the cosine measure relative to the span of the set is provided and proven to return the correct results.  In Section \ref{sec:conclusion}, the main results of this paper are summarized.

\section{Notation and Preliminaries}\label{sec:prel}
The zero vector in $\R^n$ is denoted by $\zero_n$ and the vector of all ones in $\R^n$ is denoted by $\one_n.$ When the dimension of the vector is clear, we may omit the subscript.  The $i^{th}$ coordinate vector in $\R^n$ is denoted $e_i$. 

We denote the closure of a set by $\cl$.  
We denote by $B_n(x^0;\Delta)$ the \emph{open ball centered about $x^0 \in \R^n$ with radius $0<\Delta<\infty$} and by $\overline{B}_n(x^0;\Delta)$ the \emph{closed ball centered about $x^0$ with radius $\Delta$.} That is
\begin{align*}
     B_n(x^0;\Delta)&=\left \{ x \in \R^n: \Vert x-x^0 \Vert < \Delta \right \}
     \ \mbox{ and } \\ \overline{B}_n(x^0;\Delta)&=\cl\left(B_n(x^0;\Delta)\right) = \left \{ x \in \R^n: \Vert x-x^0 \Vert \leq \Delta \right \}.
\end{align*}

Given a set of vectors $D \subseteq \R^n$ (possibly infinite), the cardinality of $D$ is denoted by $|D|$. 
The radius of the smallest ball centered at the origin containing the set $D$ is denoted by $\Delta_D$ and given by 
 \begin{align} \label{eq:radiusD}
\Delta_D&=\sup_{d \in D} \Vert d \Vert.
 \end{align}
When useful, in this paper, a set of vectors is represented as a matrix where each column represents a vector in the set.  That is
    $$D=\{d_1, d_2, \ldots, d_m\} \mbox{ is interchangable with } D = \begin{bmatrix} d_1 & d_2 & \ldots & d_m \end{bmatrix}.$$

%

The {\em dimension} of a linear subspace $L \subseteq \R^n$ is denoted by $\dim(L)$. A \emph{trivial} subspace $L\subset \R^n$ is a set of the form $L=\{\zero_n\}.$ Similarly, a set $D$ in $\R^n$ is said to be trivial if  $D=\{\zero_n\}$ or $D=\emptyset.$ 
In this paper, the linear subspace  considered are assumed to be nontrivial.

We now focus on the key definitions studied in this paper. In the remainder of this paper, the word linear may be omitted when discussing the notion of span, subspace and basis. All these notions are defined for the linear case.  

\begin{definition}[span, positive span]  Let $D \subseteq \R^n$. 
\begin{enumerate}[i.]
    \item The  {\em span} of $D$ is denoted by $\spann(D)$ and defined by
    $$\spann(D) = \left\{x \in \R^n: x =\sum_{i=1}^k \alpha_i d_i, k \in \N, d_i \in D, \alpha_i \in \R\right\}.$$
    \item The {\em positive span} of $D$ is denoted by $\pspan(D)$ and defined by
    $$\pspan(D) = \left\{x \in \R^n: x =\sum_{i=1}^k \lambda_i d_i, k \in \N, d_i \in D, \lambda_i \geq 0\right\}.$$
\end{enumerate}
\end{definition}

Note that $\spann(D)$ is always a linear subspace of $\R^n.$  Moreover, $L$ is a linear subspace of $\R^n$ if and only if $L = \spann(L)$.

The \emph{projection}  of a vector $v \in \R^n$ onto a linear subspace $L$ will be denoted by $\Proj_L v$.  We will be particularly interested in the projection onto $\spann(D)$.  For ease of writing, we provide this with the special notation $\ProjD = \Proj_{\spann(D)}$.  In Section \ref{sec:algo}, we shall make use of the well-known formula
\begin{equation*}
\ProjD v=DD^\dagger  v, 
\end{equation*}
where $D^\dagger$ denotes the \emph{Moore--Penrose pseudoinverse of $D$} \cite[Chapter 7]{roman2007}.

\begin{definition}[spanning, positive spanning]
Let $D \subseteq \R^n$ and $L$ be a subspace of $\R^n$.   
\begin{enumerate}[i.]
    \item The set $D$ is said to be a \emph{spanning set of $L,$} or \emph{spans $L$,} if and only if $\spann(D)=L.$
    \item The set $D$ is said to be a \emph{positive spanning set of $L,$} or \emph{positively spans $L$,} if and only if $\pspan(D)=L.$
\end{enumerate}\end{definition}

It is easy to prove that, if $D$ (positively) spans $L$, then $D$ must be a subset of $L$.  Indeed, one always has
    $$D \subseteq \pspan(D) \subseteq \spann(D).$$

The term (positively) independent is used to mean that removing any vector from the set changes the (positive) spanning properties of the set.

\begin{definition}[independence, positive independence]
Let $D \subseteq \R^n$.  
\begin{enumerate}[i.]
    \item The set $D$ is said to be \emph{dependent} if and only if there exists a vector $d \in D$ such that $d \in \spann(D\setminus \{d\})$. Conversely, $D$ is said to be \emph{independent} if and only if $\spann ( D\setminus{\{d\}})\neq \spann(D)$ for any $d \in D$.
    \item The set $D$ is said to be \emph{positively dependent} if and only if there exists a vector $d \in D$ such that $d \in \pspan(D\setminus \{d\})$. Conversely, $D$ is said to be \emph{positively independent} if and only if $\pspan ( D\setminus{\{d\}})\neq \pspan(D)$ for any $d \in D$.
\end{enumerate}
\end{definition}

Finally, we introduce the notion of a (positive) basis of a subspace.

\begin{definition}[basis, positive basis]
Let $D \subseteq \R^n$ and $L$ be a  subspace of $\R^n$. 
\begin{enumerate}[i.]
    \item The set $D$ is a \emph{basis of $L$} if and only if  $D$ is independent with $\spann(D)=L.$
    \item The set $D$ is a \emph{positive basis of $L$} if and only if  $D$ is positively independent with $\pspan(D)=L.$
\end{enumerate}\end{definition}

It is well known that if $D$ is a basis of $L$ then $|D| = \dim(L)$.  It is shown in \cite{Davis1954,regis2016} that the minimal cardinality of a positive basis of  a  subspace $L$ is $\dim(L)+1$ and the maximal cardinality is $2 \dim(L).$   Positive bases of these sizes are called \emph{minimal positive bases} and \emph{maximal positive bases} (respectively).  Positive bases  with cardinality strictly between $\dim(L)+1$ and $2\dim(L)$ are called  \emph{intermediate positive bases}.




We end this section by recalling the definition of  cosine measure and  defining the cosine measure relative to a subspace. The cosine measure is valuable tool to quantify how well a set covers the space $\R^n.$ We will see that the cosine measure relative to a subspace extends the value of this tool to work with a subspace of $\R^n$. We begin with the definition of cosine measure (as given in \cite{hare2023nicely,regis2021}) and the corresponding cosine vector set (from \cite{hare2020deterministic}).

\begin{definition}[cosine measure] \label{def:cosinemeasure}  Let $D \subseteq \R^n$ be a nonempty finite set of nonzero vectors. 
The \emph{cosine measure} of $D$ is defined by 
    $$\cm{D} \ =\ \min_{\substack{u \, \in \, \R^n\\\Vert u \Vert=1}} \max_{d \, \in \, D} \frac{u^\top d}{\Vert d \Vert},$$
and the  \emph{cosine vector set of $D$}, denoted by $\V(D)$, is defined by $$\V(D) \ =\ \argmin{\substack{u \, \in \, \R^n\\\Vert u \Vert=1}}\max_{d\, \in\, D} \frac{u^\top d}{\Vert d \Vert}.$$
\end{definition}

Note that the original definition of the cosine measure requires $D$ to be finite.  Another limitation of the cosine measure is that is it focuses on how well the set covers the entire space $\R^n$.  As a result, if $D$ is a positive spanning set of a proper subspace, then the cosine measure will always return $\cm{D} = 0$ (Corollary \ref{cor:cmstrictsubspacepbasis} herein).  To address these limitations, we introduce the {\em cosine measure relative to a subspace}, which provides information on how well a set covers a linear subspace.  Simultaneously, we define a corresponding cosine vector set and allow for both concepts to be well-defined for infinite sets.

However, since an empty set, or the zero vector, provides no coverage of $\R^n$, the cosine measure relative to a subspace is defined for nonempty sets of nonzero vectors. If a set is empty, the cosine measure (relative to a subspace) should be assumed to be undefined.  If a set contains the zero vector, then one should remove the vector and work with the remaining set. 
If the remaining set is empty, then the result would again be undefined.
\begin{center}\fbox{%
    \parbox{0.8\textwidth}{%
        For the remainder of this paper we shall assume that the set $D$ is a nonempty set of nonzero vectors.
    }%
}\end{center}

\begin{definition}[cosine measure relative to $L$] \label{def:cmRestricted}
Let $L \subseteq \R^n$ be a nontrivial linear subspace and let
    $D \subseteq \R^n$ be a nonempty set of nonzero vectors.
The \emph{cosine measure of  $D$ relative to $L$} is denoted by $\CM_L(D)$ and defined by 
\begin{align} \label{eq:cmsubspace}
    \CM_{L}(D)\ = \ \min_{\substack{u \, \in \, L\\\Vert u \Vert=1}} \sup_{\substack{d \, \in \, D}} \frac{u^\top d}{\Vert d \Vert},
\end{align}
and the \emph{cosine vector set of $D$ relative to $L$}, denoted by $\V_{L}(D)$, is defined by $$\V_{L}(D)\ =\ \argmin{\substack{u \, \in \, L\\\Vert u \Vert=1}}\sup_{\substack{d\, \in\, D}} \frac{u^\top d}{\Vert d \Vert}.$$
\end{definition}

Notice that the cosine measure is a special case of the cosine measure relative to $L$.  Indeed, if $L=\R^n$ and $|D|$ is finite, then the cosine measure relative to $L$ returns the classical cosine measure: $\CM_{\R^n}(\cdot) = \cm{\cdot}$.  
Indeed, if $D$ is finite, then the $\sup$ can be replaced by $\max$ in \eqref{eq:cmsubspace}.  
If $D$ is a trivial set or if $L$ is a trivial subspace of $\R^n,$ then the cosine measure and cosine vector set (relative to $L$) are undefined.  
When  $D$ is a nonempty set of nonzero vectors  and $L$ is a nontrivial subspace, the cosine measure relative to a subspace  returns a value in $[-1,1]$. Besides,
the cosine vector set of $D$ relative to $L$ is well-defined and nonempty.  Notice that when $D$ is an infinite set, then the cosine vector set may be infinite. For example, consider the infinite set $D=\{x \in \R^2: \Vert x \Vert=1\}.$ Then $\V_{\R^2}(D)=D.$ It is also possible to obtain an infinite cosine vector set  when  $D$ is a finite set (for instance, see $\V(D_2)$ in Example \ref{ex:valueCMrelative}).
\begin{remark}
Given a nonempty set of nonzero vectors $D$ and a nontrivial linear subspace $L$ in $\R^n$ such that  $D \not\subseteq L,$ one might consider projecting each normalized vector $d/\Vert d \Vert \in D$ onto $L$ before computing the cosine measure of $D$ relative to $L.$ 
Denote by $\widetilde{D}$ the set  obtained  from these projections and by $\widetilde{d}$ a vector in the set $\widetilde{D}$. Note that the set $\widetilde{D}$ might contain the zero vector. If the computations are done using the projected set $\widetilde{D},$  then it is crucial  to keep the   zero vector in the set $\widetilde{D}.$ The cosine measure of $D$ relative to $L$ can be  computed as follows whenever we work from the projected set $\widetilde{D}$ rather than using the original set $D$: 

\begin{align}
\CM_L(D) &= \min_{\substack{u \, \in \, L\\ \Vert u \Vert=1}} \sup_{d \, \in \, D} \frac{u^\top d}{\Vert d \Vert } \notag \\
&= \min_{\substack{u \, \in \, L\\ \Vert u \Vert=1}} \sup_{d \, \in \, D} u^\top \left (\proj_L\frac{d}{\Vert d \Vert } + \left  (\frac{d}{\Vert d \Vert} - \proj_L \frac{d}{\Vert d \Vert}  \right )\right ) \notag \\ 
&= \min_{\substack{u \, \in \, L\\ \Vert u \Vert=1}} \sup_{\widetilde{d} \, \in \, \widetilde{D}} u^\top \widetilde{d}, \label{eq:cmProjectedalternate} 
\end{align}
 as $u^\top  \left (\frac{d}{\Vert d \Vert} - \proj_L\frac{d}{\Vert d \Vert }\right ) = 0$ for any $u \in L$.
\end{remark}

The following example shows the importance of keeping the zero vector in the projected set $\widetilde{D}$  if the computations are done  from the projected set $\widetilde{D}.$

\begin{example}
Let $D=\{e_1,e_2\} \subseteq \R^2$ and $L=\{x \in \R^2:x_2=0\}.$ Then $$\CM_L(D) \ = \ 0.$$ The projection of $e_1$ onto $L$ is equal to $e_1$ and the projection of $e_2$ onto $L$ is equal to $\zero.$ We obtain the projected set $\widetilde{D}=\{e_1, \zero\}.$ Computing the cosine measure of $D$ relative $L$  as defined in Equation \eqref{eq:cmProjectedalternate}, we obtain $\CM_L(D)=0.$  However, if one removes the zero vector from $\widetilde{D},$ the final result would be 
    $$\min_{\substack{u \, \in \, L\\ \Vert u \Vert=1}} \sup_{\widetilde{d} \, \in \, \widetilde{D}\setminus\{\zero\}} u^\top \widetilde{d} \ = \ -1 \ \ne \ \CM_L(D).$$
\end{example}

We shall often be interested in the special case when $L = \spann(D)$.  For ease of writing, when the linear subspace considered is the span of a set,  we only write the set name and omit the word span  in the subscript. More precisely,
    $$\CM_{D}(\cdot)= \CM_{\spann(D)}(\cdot) \quad \mbox{and} \quad \CV_{D}(\cdot) = \V_{\spann(D)}(\cdot).$$
Finally, some analysis will require examining the active set of the cosine measure relative to $L$.

\begin{definition}[active set]\label{def:allactivity} 
Let $L \subseteq \R^n$ be a nontrivial linear subspace and 
    $D \subseteq \R^n$ be a nonempty set of nonzero vectors.

The \emph{active set for $D$ at a cosine vector $u \in \V_{L}(D)$  relative to $L$} is denoted by  $\A_{L}(D,u)$ and defined by $$\A_{L}(D,u)\ =\ \left \{d \in D\,:\,\frac{d^\top u}{\Vert d \Vert}=\CM_{L}(D) \right \}.$$
The \emph{active set for $D$ relative to $L$} is denoted by $\AL_{L}(D)$, and  defined by $$\AL_{L}(D)\ =\ \bigcup_{u \, \in \, \V_L(D)} \A_{L}(D,u).$$
\end{definition}

As above, we are mostly interested in the special case of $L = \spann(D)$.  For this case we use the special notation
    $$\A_D(\cdot, \cdot) = \A_{\spann(D)}(\cdot, \cdot).$$
In Section \ref{sec:algo}, we will see that $\A_D(D,u)$ contains a basis of $\spann(D)$ whenever $\CM_D(D)>0$ and $D$ is a nonempty closed set of nonzero vectors.  If $D$ is not closed, then this result may not hold.  For example,  consider the set $D=S_2 \setminus \{x \in \R^2: x_1^2+x_2^2=1, x_1\geq 0,x_2 \geq 0\}.$ Then $\CM_D(D)=1/\sqrt{2}$ and  the cosine vector is $u=\begin{bmatrix} 1/\sqrt{2}&1/\sqrt{2}\end{bmatrix}^\top.$  The only two vectors returning the cosine measure value are $v=\begin{bmatrix} 0&1 \end{bmatrix}^\top$ and $w=\begin{bmatrix} 1&0\end{bmatrix}^\top$.  Neither are in $D$, so $\A_D(D,u)=\emptyset.$  If we change $D$ to $S_2\setminus \{x \in \R^2: x_1^2+x_2^2=1, x_1>0,x_2\geq 0\}$, then  $\A_D(D,u)=\left\{[1, 0]^\top\right\}$, so we can have $\A_D(D,u)$ nonempty, but not containing a basis.  Finally notice that if we take the closure of either of the above sets, we have $\CM_{\cl(D)}(\cl(D))=1/\sqrt{2}$ and  $\A_{\cl(D)}(\cl(D),u)=\left\{[1, 0]^\top, [0,1]^\top\right\}$.

To avoid this situation, we sometimes assume that the set considered is closed in the main results of this paper (note that a finite set of vectors is necessarily a closed set). This assumption does not represent a strong assumption   since taking the closure of a set $D$  does not change the value of the cosine measure.

\section{Properties of a positive spanning set of a linear subspace} \label{sec:properties}

In this section, properties of positive spanning sets of a linear subspace are investigated.  The results are developed to consider  both the cases where $|D|$ is finite and $|D|$ is infinite.   We begin by showing  that a set $D$ contained in a linear  subspace $L$ is a positive spanning set of $L$ if and only if $\pspan(D)=\spann(D)=L.$

 \begin{lemma} \label{lem:LisSD} Let $D \subseteq L$ where $L$ is a nontrivial linear subspace of $\R^n$.  Then $D$ positively spans $L$ if and only if $\pspan(D)=\spann(D)=L.$ 
 \end{lemma}
\begin{proof} ($\Rightarrow$) 
Suppose that $D$ positively spans $L$; i.e., $\pspan(D)=L$.  This implies $D\subseteq L$ and therefore $\pspan(D) \subseteq \spann(D) \subseteq L = \pspan(D)$.
\\
$(\Leftarrow)$ 
Conversely, if $\pspan(D)=\spann(D)=L$, then $D$ positively spans $L$ by definition.\qed
\end{proof}

We next show that if $|D|$ is infinity and the positive span is a linear subspace, then the positive span can be created through a finite subset of $D$.  We use the following lemma. 

\begin{lemma}\cite[Theorem 2.11]{Brown1988} \label{lem:containsBasisSD}Let $D \subseteq \R^n$ be a nonempty (possibly infinite) set of nonzero vectors.  Then $D$ contains a basis of $\spann(D)$.
\end{lemma}

\begin{corollary} \label{cor:containFinitePSS}Let $D \subseteq \R^n$ be a nonempty (possibly infinite) set of nonzero vectors and $L$ be a linear subspace of $\R^n$.  If $\pspan(D)=L$, then there exists a finite subset of vectors $C \subseteq D$ such that $\pspan(C)=L$.
\end{corollary}
\begin{proof} 
If $D$ is a finite set, then take $C=D.$  

Suppose that $D$ is an infinite set and $\dim(L)\geq 1$.  By Lemma \ref{lem:LisSD}, $\pspan(D) = L$ implies $L=\spann(D)$.  By Lemma \ref{lem:containsBasisSD}, the set $D$ contains a basis of $\spann(D)$.  Denote this basis by $B=\left \{ b_1, b_2, \ldots, b_m \right \}$ where $m = \dim(\spann(D))$.  Define $v=-\sum_{i=1}^m b_i.$  Since $v \in \spann(D)=\pspan(D)$,  by Carathéodory's Theorem \cite[Theorem 1.3.6]{hiriart2004fundamentals}, the vector $v$ can be written as $$v=\alpha_1 d_1 + d_2 +\dots +\alpha_k d_k$$ where $d_i \in D, \alpha_i\geq 0,$ and $k$ is an integer greater or equal than 1.  Since $v  \in \pspan(B \cup \{v\})=\spann(D)$ \cite[Theorem 5.1]{regis2016}, we have $$\pspan(B \cup \left \{d_1, d_2 \ldots, d_k  \right \})=\pspan(B \cup \{v\})=\spann(D)=L.$$
The set $B \cup \left \{d_1, d_2 \ldots, d_k  \right \}$ is finite, so the proof is complete. \qed
\end{proof}

Note that the above result requires $D$ to be a positive spanning set.  Indeed, consider the {\em floating ring set} defined by 
$D=\{ x \in \R^3: (x_1)^2+(x_2)^2=1, x_3=1\}.$
It can be shown that the set $D$ is positively independent \cite{hare2016cardinality}.  As such, for any finite subset $C \subseteq D,$ we have $\pspan(C) \neq \pspan(D).$

This inspires the following result about positive independence for an infinite set of vectors.

\begin{corollary}\label{cor:infiniteSetDoesNotPspanL}
Let $D \subseteq \R^n$ be an infinite set of vectors. If $D$ is positively independent, then  $\pspan(D) \neq \spann(D).$ 
\end{corollary}
\begin{proof} Suppose that $D$ is positively independent and 
 $\pspan(D)=\spann(D).$ By Corollary \ref{cor:containFinitePSS}, there exists a finite subset subset of $D$, say $C$,  such that $\pspan(C)=\spann(D)$.  Since $D$ is positively independent, we have $\pspan(C) \neq   \pspan(D)$ which  creates a contradiction.\qed
\end{proof}

Notice that Corollary \ref{cor:infiniteSetDoesNotPspanL} is specific to infinite sets.  Indeed, a finite set of positively independent vectors $D$ such that $\pspan(D) = \spann(D)$ is a positive basis. 

We now extend \cite[Theorem 2.5]{regis2016} to the cases of infinite sets.

\begin{theorem} \label{thm:equivalencepss}Let $D \subseteq \R^n$ be a nonempty (possibly infinite) set of nonzero vectors.  Then the following are equivalent.
\begin{enumerate}[(i)]
    \item The set $D$ is a positive spanning set of $\spann(D).$ 
    \item For any $d \in D$, $-d$ is in $\pspan(D\setminus \{d\})$.
\end{enumerate}
Moreover, these imply that there exists a finite subset $C\subseteq D$ such that $\pspan(C)=\pspan(D)$.

Suppose that $C$ is a finite subset of $D$ such that $\pspan(C)=\pspan(D)$ (if $D$ is finite, then we may take $C=D$).  Let $s = |C|\geq 2$.  Then the following are equivalent.
\begin{enumerate}[(i)]\setcounter{enumi}{2}
    \item The set $D$ is a positive spanning set of $\spann(D).$ 
    \item The set $C$ is a positive spanning set of $\spann(D).$ 
    \item There exists $\alpha \in \R^s$, such that $\alpha>0$ and  $C\alpha = \zero_n$.
    \item There exists $\beta \in \R^s$, such that $\beta \geq 1$ and $C\beta = \zero_n$.
    \item There exists $\gamma \in \R^s$ such that $\gamma \geq 0$ and $C\gamma=-C\one_s$.
\end{enumerate}
where items (v), (vi), and (vii), interpret $C$ as a matrix in $\R^{n \times s}$.
\end{theorem}

\begin{proof} Parts (i) and (ii) are equivalent by \cite[Theorem 2.5]{regis2016}.  
Suppose that $D$ is a positive spanning set of $\spann(D)$.
By Corollary \ref{cor:containFinitePSS}, part (i) implies that there exists a finite subset $C\subseteq D$ such that $\pspan(C)=\spann(D).$  The proof of equivalence of parts (iii) to (vii) is identical to the proof for finite sets provided in \cite[Theorem 2.5]{regis2016}. \qed 
\end{proof}

Our next example shows the importance of $\pspan(C)=\pspan(D)$ in the second half of Theorem \ref{thm:equivalencepss}.   

\begin{example} Consider $D_{\infty}=\{ d : \|d\|=1, d_1 \geq 0\} \subseteq \R^2$.  The set $C = \{e_2, -e_2\}$ is contained in $D_\infty$, but $\pspan(C) \neq \pspan(D_\infty)$.  Notice that $C$ satisfies parts (v), (vi), and (vii) of Theorem  \ref{thm:equivalencepss}, but $C$ is not a positive spanning set of $\spann(D_\infty)$. 
\end{example}


In practice, Theorem \ref{thm:equivalencepss}(vii) is useful to decide if a given set is a positive spanning set.  Theorem \ref{thm:equivalencepss}(vii) can also be used to show how to extend a set that is not a positive spanning set into a positive spanning set by adding only one vector. In the finite case, these results   have  been known for a long time and are discussed in \cite{Davis1954}. Hence, the value of the following proposition is to show that it also holds for infinite sets.

\begin{proposition} \label{prop:extendPSSSD}Let $D \subseteq \R^n$ be a nonempty (possibly infinite) set of nonzero vectors.   Let $m = \dim(\spann(D))$ and $B=\{b_1, \cdots, b_m\} \subseteq D$ be a basis of $\spann(D)$.  Define the vector $w=-\sum_{j=1}^m b_j.$ Then
\begin{enumerate}[(i)]
    \item the set $D'=D \cup \{w\}$ is a positive spanning set of $\spann(D)$, and 
    \item the set $D''=D \cup -B$ is a positive spanning set of $\spann(D)$.
\end{enumerate} 
\end{proposition}

\begin{proof}
If $\pspan(D)=\spann(D),$ then the result holds immediately. So suppose that $\pspan(D) \neq \spann(D).$ 

\noindent (i) Let $C = B \cup \{w\}$.  By \cite[Thm 6.4]{audet2017derivative}, $\pspan(C)=\spann(B)=\spann(D)$.  We also have that the sum of the vectors of $C$ is the null vector.
Therefore, Theorem \ref{thm:equivalencepss}(v) ensures that $\pspan(D')=\spann(D)$. 
\\
(ii) The result follows similarly using $C = B \cup (-B)$. 
 \qed
\end{proof}

One of the most useful properties of a positive spanning set of $\R^n$ is  that  for any nonzero vector $v \in \R^n,$ there exists a vector $d$ in the positive spanning set for which 
\begin{align*}
v^\top d&>0. 
\end{align*}
It follows that given $f \in \mathcal{C}^1$, a positive spanning set of $\R^n$ contains a descent direction of $f$ at any point $x^0 \in \R^n$ where $\nabla f(x^0) \neq \zero_n$ \cite{conn2009introduction}. This result can be generalized to the linear subspaces and infinite sets.  To do this, we apply the following lemma.

\begin{lemma}\label{lem:rockafellar}Let $D \subseteq \R^n$ be a nonempty (possibly infinite) set of nonzero vectors contained in a linear subspace $L$.  If $\pspan(D) \neq L,$ then all the vectors in $D$ are contained in  a closed half-space of $L.$ 
\end{lemma}

\begin{proof}Follows immediately from  \cite[Corollary 11.7.3]{rockafellar1970}. \qed
\end{proof}

\begin{proposition} \label{prop:dotProdPos} Let $D \subseteq \R^n$ be a nonempty (possibly infinite) set of nonzero vectors contained in a linear subspace $L$.
Then $\pspan(D)=L$ if and only if for any nonzero vector $v \in L$ there exists a vector $d \in D$ such that 
\begin{equation}
v^\top d>0. \label{eq:dotProdPos}
\end{equation}
\end{proposition}
\begin{proof}
($\Rightarrow$) Suppose that $\pspan(D)=L$.  Let $v$ be a nonzero vector in $L$.  As $v \in \pspan(D)$, $v$ can be written as
$$v=\alpha_1 d_1+\alpha_2d_2+\dots+\alpha_m d_m,$$
where $\alpha_j\geq 0$ for all $j \in \{1, 2, \ldots, m\}.$
It follows that
$$0<v^\top v=\alpha_1 d_1^\top v+\alpha_2 d_2^\top v+ \dots +\alpha_m d_m^\top v.$$
Hence, we must have $d_j^\top v>0$ for at least one $j \in \{1, 2, \ldots, m\}.$
\\
($\Leftarrow$) 
Conversely, suppose that given any $v \in L \setminus \{\zero_n\}$ there exists a $d \in D$ such that $v^\top d>0$. 
If $\pspan(D) \neq L$, then by Lemma \ref{lem:rockafellar} there exists $h \in \spann(D) \setminus\{\zero_n\}$ such that all vectors in $D$ are contained in $\{w \in \spann(D) : w^\top h \leq 0\}$.  This leads to a contradiction  by taking $v=h$.
\qed
\end{proof}

In Proposition \ref{prop:dotProdPos}, if we  let the nonzero vector  $v$ be in $\R^n$ rather than $L,$ then the result  does not necessarily hold. The forward direction is true if we replace the strict inequality in \eqref{eq:dotProdPos} by an inequality. Proposition \ref{prop:dotProdNonneg} provides a proof of this claim. 

\begin{proposition} \label{prop:dotProdNonneg}Let $D \subseteq \R^n$ be a nonempty (possibly infinite) set of nonzero vectors contained in a linear subspace $L$. If $\pspan(D)=L$, then for any vector $v \in \R^n,$ there exists a vector $d \in D$ such that 
 \begin{equation*}
 v^\top d \geq 0. \label{eq:dotProdNonneg}
 \end{equation*}
\end{proposition}
\begin{proof}
Let $v \in \R^n.$ An orthogonal decomposition of $v$ is given by $v=\ProjD v+v^\perp$ where $v^\perp$ is in the subspace  of all vectors in $\R^n$ that are orthogonal to $L.$ Note that $v^\top d = (\ProjD v+ v^\perp)^\top d= (\ProjD v)^\top d$ for all $d \in D$.
If $\ProjD v=\zero_n$, then $v^\top d=0$ for all $d \in D$.  If $\ProjD  v \neq \zero_n,$ then by Proposition \ref{prop:dotProdPos} there exists $d \in D$ such that $(\ProjD  v)^\top d> 0.$ \qed
\end{proof}

The converse of Proposition \ref{prop:dotProdNonneg} is not necessarily true. For example, consider the set $D=\{e_1, -e_1, e_2\} \subseteq \R^2$.  For this set, for any vector $v \in \R^2$, there exists $d \in D$ such that $v^\top d \geq 0$, however $\pspan(D) \neq \spann(D)$.

In the next section,  the notion of cosine measure relative to a linear subspace is investigated. 

\section{Cosine measure relative to a subspace}\label{sec:CMrelativetoL}

Definition \ref{def:cmRestricted} extends the definition of the cosine measure to include the idea of a cosine measure relative to a subspace.  In this section we explore the basic properties of this new definition.  As we are most interested in the cosine measure relative to $\spann(D)$, many results focus on $\CMD$.

Since the original definition of the cosine measure assumes $D$ is finite (Definition \ref{def:cosinemeasure}), results involving the cosine measure include the assumptions that $D$ is finite.

We begin with the obvious relationship between the cosine measure and the cosine measure relative to a subspace. 

\begin{proposition} \label{prop:inequalityCmAndCmsd}
Let $D \subseteq \R^n$ be a nonempty finite set of nonzero vectors, $L$ be a nontrivial linear subspace of $\R^n$.   
Then $$\CM(D) \leq   \CM_{L}(D)$$
with equality if and only if $\V(D) \cap L \neq \emptyset$.
\end{proposition}
\begin{proof}
Since  $L$ is a nontrivial subspace  of $\R^n$, the definitions of the cosine measures ensure that 
    \begin{align*}
    \CM(D)\ =\ 
        \min_{\substack{u \, \in \, \R^n\\\Vert u \Vert=1}} \max_{\substack{d \, \in \, D }} \frac{u^\top d}{\Vert d \Vert} 
        \ \leq \ \min_{\substack{u \, \in \, L\\ \Vert u \Vert=1}} \max_{\substack{d \, \in \, D }} \frac{u^\top d}{\Vert d \Vert}
        \ =\ \CM_{L}(D).
    \end{align*}
    Now, suppose that $\CM(D)=\CM_{L}(D).$ Let $u^* \in \V_{L}(D).$ Then $u^*$ must be in $\V(D)$, so $\V(D) \cap L \neq \emptyset.$ Conversely, suppose that $\V(D) \cap L \neq \emptyset.$   Let $u^\# \in \V(D) \cap L.$ Then 
    \begin{align*}
        \CM(D) 
        \ =\ \max_{\substack{d \, \in \, D }} \frac{(u^\#)^\top d}{\Vert d \Vert}
        \geq \  \min_{\substack{u \, \in \, L \\ \Vert u \Vert=1}} \max_{\substack{d \, \in \, D }} \frac{u^\top d}{\Vert d \Vert} \ =\ \CM_{L}(D). 
    \end{align*} Since $\CM(D) \leq \CM_L(D)$, we must have $\CM(D)=\CM_{L}(D).$ \qed
\end{proof}

\subsection{Relating $\pspan(D)$ and $\CMD$}\label{ssec:}

We now turn our attention to how the cosine measure relative to a subspace provides knowledge about positive spanning properties.  We begin with an example computing the cosine measure relative to a subspace.  We will return to this example after each result to illustrate what the information the cosine measure relative to a subspace provides. In the following example, it is relatively easy to find the  exact value of the cosine measure relative to a subspace since the subspaces considered are either  one-dimensional or two-dimensional.  A deterministic algorithm will be provided in Section \ref{sec:algo}.

\begin{example}\label{ex:valueCMrelative}
Consider the sets of directions 
    $D_1=\{e_1, -e_1, e_2, -e_2\} \subseteq \R^3$ and 
    $D_2=\{e_1, -e_1, e_2\} \subseteq \R^3$
    and the subspaces 
        $L = \{x \in \R^3 : 4x_1 = 3x_2, x_3 = 0\}$ and 
        $M=\{x \in \R^3 : x_1 = x_3=0\}$.  
    First, note that 
   $$\begin{array}{rll}
       &\qquad \cm{D_1}=0, &\qquad \V(D_1) = \{e_3, -e_3\} \\ 
       \text{and}  & \qquad
        \cm{D_2}=0, &\qquad \V(D_2) = \{v = (0, v_2, v_3) \, :\, \|v\| =1, \, v_2 \leq 0\}. 
   \end{array}$$

The shaded regions in Figure~\ref{fig-ex} represent the positive span of $D_1$ on the left and of $D_2$ on the right.
Both figures also show the subspaces $L$ and $M$.

\begin{figure}[htb!]
\begin{minipage}{0.5\textwidth}
\begin{tikzpicture}[scale=1]
  \begin{axis}[
	width=0.8\textwidth,
	height=0.8\textwidth,
	axis lines = none,
    xmin = -10, xmax = 10,
    ymin = -5, ymax = 5.3,
]
	\draw[ultra thin, fill=gris9]  (-10,-2) -- (-4,2) -- (10,2) -- (4,-2) -- cycle;
	\draw[ultra thick, gris4] (-7,0) -- (7,0); 	\node[gris4] at (-5.5,0.2) {$^M$}; 
	\draw[ultra thick, gris6] (9.2,1.5) -- (-9.2,-1.5); 	\node[gris6] at (-8.,-1.8) {$^L$}; 
	\node[gris7] at (5.7,-2.7) {$^{\pspan(D_1)}$}; 
	\draw[thin, ->] (3,2) -- (-3.6,-2.4); 	\node at (-3.6,-3.) {$-e_1$};
 	\draw[thin, ->] (-3.6,-2.4) -- (3.6,2.4); 	\node at (4.2,2.8) {$e_1$};
	\draw[thin, ->] (-7,0) -- (8,0); 	 \node at (7.6,-0.6) {$e_2$};
	\draw[thin, ->] (8,0) -- (-8,0); 	 \node at (-9.4,0.4) {$-e_2$};
	\draw[thin, <->] (0,-3) -- (0,4.2); 
\end{axis};
\end{tikzpicture}
\end{minipage}
\hspace*{-20mm}
\begin{minipage}{0.5\textwidth}
\begin{tikzpicture}[scale=1]
  \begin{axis}[
	width=0.8\textwidth,
	height=0.8\textwidth,
	axis lines = none,
    xmin = -10, xmax = 10,
    ymin = -5, ymax = 5.3,
]
	\draw[ultra thin, fill=gris9]  (-3,-2) -- (3,2) -- (10,2) -- (4,-2) -- cycle;
	\draw[ultra thick, gris4] (-7,0) -- (7,0); 	\node[gris4] at  (-5.5,0.2) {$^M$}; 
	\draw[ultra thick, gris6] (9.2,1.5) -- (-9.2,-1.5); 	\node[gris6] at (-8.,-1.8) {$^L$}; 
	\node[gris7] at (5.7,-2.7) {$^{\pspan(D_2)}$}; 
	\draw[thin, ->] (3,2) -- (-3.6,-2.4); 	\node at (-3.6,-3.) {$-e_1$};
 	\draw[thin, ->] (-3.6,-2.4) -- (3.6,2.4); 	\node at (4.2,2.8) {$e_1$};
	\draw[thin, ->] (-7,0) -- (8,0); 	 \node at (7.6,-0.6) {$e_2$};
	\draw[thin, <->] (0,-3) -- (0,4.2); 
\end{axis};
\end{tikzpicture}
\end{minipage}
\caption{Illustrations of the positive spans of $D_1$ and $D_2$}
\label{fig-ex}
\end{figure}
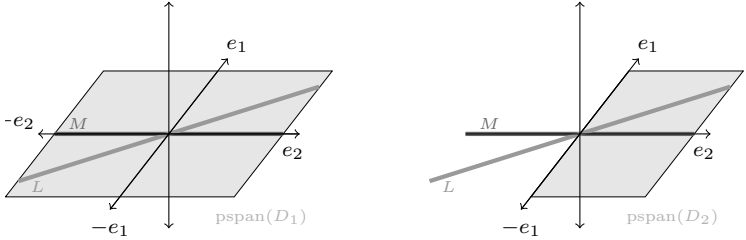

    There are exactly two unit vectors in 
        $L$:        $$u^1 = \begin{bmatrix}
            0.6 \\ 0.8 \\ 0
        \end{bmatrix} \quad \mbox{and} \quad 
        u^2 = \begin{bmatrix}
            -0.6 \\ -0.8 \\ 0
        \end{bmatrix}.$$
    The cosine measures and cosine vectors relative to the subspace $L$ are
   $$\begin{array}{rll}
        &\qquad \CM_{L}(D_1)\ =\ \min\left\{ \displaystyle\max_{d \in D_1}  u^\top d \, :\, 
            u \in \{u^1,u^2\} \right\} \ =\ 0.8,
            & \qquad \V_{L}(D_1) \ = \ \{ u^1, u^2 \} \\
       \text{and} & \qquad
        \CM_{L}(D_2) \ = \ \min\left\{ \displaystyle\max_{d \in D_2} u^\top d\, :\, 
            u \in \{u^1,u^2\} \right\} \ =\ 0.6,
            & \qquad \V_{L}(D_2) \ = \ \{ u^2 \}.
   \end{array}$$
    In the subspace $M,$ there are exactly two unit vectors: $e_2$ and $-e_2.$ We obtain
   $$\begin{array}{rll}
        &\qquad \CM_{M}(D_1)\ =\ \min\left\{ \displaystyle \max_{d \in D_1} u^\top d\, :\, 
            u \in \{e_2, -e_2 \} \right\} \ =\ 1,
                    & \qquad \V_{M}(D_1) \ = \ \{ e_2, -e_2 \} \\
       \text{and} & \qquad
        \CM_{M}(D_2) \ = \ \min\left\{ \displaystyle \max_{d \in D_2} u^\top d\, :\, 
            u \in \{e_2, -e_2 \} \right\} \ =\ 0,
            & \qquad \V_{M}(D_2) \ = \ \{ -e_2 \}.
   \end{array}$$
    Considering the subspaces $\spann(D_1)$ and $\spann(D_2),$ we find 
   $$\begin{array}{rll}
        &\qquad \CM_{D_1}(D_1)=1/\sqrt{2},
        & \qquad \V_{D_1}(D_1) \ = \ \{ (v_1, v_2, 0) \, : \, |v_1| = |v_2| = 1/\sqrt{2} \} \\
       \text{and} & \qquad
        \CM_{D_2}(D_2)=0,
                & \qquad \V_{D_2}(D_2) \ = \ \{ -e_2 \}. \\
   \end{array}$$
\end{example}
    
Given a nontrivial  finite set $D\subseteq\R^n,$ it is known that $\CM(D)>0$ if and only if $D$ positively spans $\R^n$ \cite[Theorem 4.2]{regis2021}. Note that  since $\pspan(D) \subseteq \spann(D)$,  the  previous result  can be expressed as follows: given a nontrivial finite set $D \subseteq \R^n,$   we have $\CM(D)>0$ if and only if $\pspan(D)=\spann(D).$  Thus, neither $D_1$ nor $D_2$ in Example \ref{ex:valueCMrelative} is a positive spanning set of $\R^3$.  The cosine measure provides no further information. 

The next proposition shows the ability of the cosine measure relative to a subspace to detect if the subspace is contained in the positive span of $D$. The result holds for infinite sets and hence, can be viewed as an extension of \cite[Theorem 3.1]{Davis1954}. 

\begin{proposition} \label{prop:pssSDiffcmGreaterEqualZero}
Let $D \subseteq \R^n$ be a nonempty (possibly infinite) set of nonzero vectors contained in a linear subspace $L$.  Then $\pspan(D)=L$ if and only if $\CM_{L}(D)>0.$  In particular, $\pspan(D)=\spann(D)$ if and only if $\CMD>0$.\end{proposition}
\begin{proof}
($\Rightarrow$) Suppose that $\pspan(D)=L.$ By Proposition \ref{prop:dotProdPos}, For any unit vector $u \in L$, there exists $d \in D$ such that $u^\top d>0.$ Hence, $\CM_L(D)>0.$

($\Leftarrow$) Conversely, suppose that $\CM_L(D)>0. $ This means that for any unit vector $u \in L,$ there exists $d \in D$ such that $u^\top d>0.$ Hence for any nonzero vector $v$ in $L,$ we have $v^\top d=\Vert v \Vert u^\top d>0.$ The result now follows from Proposition \ref{prop:dotProdPos}. \qed
\end{proof}

Note that the previous proposition assumes the set $D$ is contained in the subspace $L.$ If it is not the case, then the result does not necessarily hold.  For instance, In Example \ref{ex:valueCMrelative} we have $\CM_L(D_2)=0.6,$ but $\pspan(D_2)$ is not equal to $L$ nor $\spann(D_2).$ Since $D_1 \subseteq \spann(D_1)$  and $D_2 \subseteq \spann(D_2),$ we  can conclude $\pspan(D_1)=\spann(D_1)$ and $\pspan(D_2) \neq \spann(D_2).$


Next we prove the claim in the introduction: ``if $D$ is a positive spanning set of a  {\em proper} subspace, then the cosine measure will always return $0$''.

\begin{corollary}\label{cor:cmstrictsubspacepbasis}
Let $D \subseteq \R^n$ be a nonempty finite set of nonzero vectors.  Suppose that $\spann(D) \neq \R^n$.  If $\pspan(D)=\spann(D)$, then $\CM(D)=0$.
\end{corollary}

\begin{proof}
Select any $v \in \R^n$ with $\|v\|=1$.  Let $v_D = \ProjD v$ and $v_{D^\perp} = v - v_D$.  Since $\pspan(D)=\spann(D)$, Proposition \ref{prop:pssSDiffcmGreaterEqualZero} implies $\CMD> 0$, which further implies 
    $$\max_{d \, \in \, D} \frac{(v_D)^\top d}{\Vert d \Vert} \ \geq \ 0,$$
with equality only if $v_D = \zero_n$.  Since $v_{D^\perp}$ is in the orthogonal subspace to $D$, 
    $(v_{D^\perp})^\top d = 0$ for all $d \in D$.
Thus, 
    $$\max_{d \, \in \, D} \frac{(v)^\top d}{\Vert d \Vert} 
    \ =\ 
    \max_{d \, \in \, D} \frac{(v_D)^\top d+(v_{D^\perp})^\top d}{\Vert d \Vert}=\max_{d \in D} \frac{(v_D)^\top d}{\Vert d \Vert}
    \ \geq \ 0,$$
with equality only if $v_D = \zero_n$.  Selecting $v$ to be in  the orthogonal subspace to $\spann(D)$ now demonstrates $\CM(D)=0$. \qed
\end{proof}

In Example \ref{ex:valueCMrelative}, $\spann(D_1) \neq \R^3$ and $\CM(D_1)=0$, so Corollary \ref{cor:cmstrictsubspacepbasis} allows for the possibility that $\pspan(D_1)=\spann(D_1)$.  However, notice that this is not sufficient to ensure that result, as $D_2$ in  Example \ref{ex:valueCMrelative} also has $\CM(D_2)=0$, but $\pspan(D_2) \neq \spann(D_2)$. Hence, Corollary \ref{cor:cmstrictsubspacepbasis} cannot be made into an `if and only if' statement.

The next theorem further investigates the relation  between the value of the cosine measure  and the value of the cosine measure relative to the subspace $\spann(D).$
 
\begin{theorem} \label{thm:cmNotEqualToCmSd}
Let $D$ be a nonempty finite set of nonzero vectors in $\R^n$. If $\CM(D) \neq \CMD,$ then $\pspan(D) = \spann(D) \neq \R^n$.
\end{theorem}
\begin{proof} If $\CM(D) \neq \CMD,$ then Proposition \ref{prop:inequalityCmAndCmsd} using $L=\R^n$ implies that $\spann(D)\neq \R^n$.

For eventual contradiction, suppose that $\CM(D) \neq \CMD$ and $\pspan(D)\neq\spann(D)$. From Proposition \ref{prop:inequalityCmAndCmsd}, $\CM(D) \neq \CMD$ implies that $\CM(D) < \CMD$.  From Proposition \ref{prop:pssSDiffcmGreaterEqualZero}, $\pspan(D)\neq\spann(D)$ implies that $\CMD\leq 0$.  Hence, we have 
\begin{align} \label{eq:ineqaulitycmandcmsd}
    \CM(D) &< \CMD \leq 0.
\end{align}
Let $v \in \V(D)$ and notice that Proposition \ref{prop:inequalityCmAndCmsd} implies $v \notin \spann(D)$.  Let $v_D = \ProjD v$ and notice that $v^\top d = v_D^\top d$ for all $d \in D$.  This implies that 
 $$\CM(D) = \max_{d \, \in D} v^\top \frac{d}{\|d\|} =  \max_{d \, \in D} v_D^\top \frac{d}{\|d\|},$$
which further implies $v_D \neq \zero_n$.  
Define $\ell= \Vert v_D \Vert$.  Considering $v_D \in \spann(D)$ and $\| \frac{1}{\ell}v_D \| = 1$, we obtain
    $$\CMD \leq \max_{d \, \in D} \left(\frac{1}{\ell} v_D\right)^\top \frac{d}{\|d\|} 
    = \frac{1}{\ell} \CM(D) \ <\  \CM(D),$$
where the last inequality comes from $0<\ell< 1$ and $\CM(D)<0$.  This contradicts inequality \eqref{eq:ineqaulitycmandcmsd}, so the initial supposition cannot hold.\qed 
\end{proof}

Returning to Example \ref{ex:valueCMrelative}, notice that $\CM(D_1)=\CM_{D_1}(D_1)$, and indeed we have $\spann(D_1)=\pspan(D_1)\neq\R^3$.  In  Example \ref{ex:valueCMrelative}, we find the $\CM(D_2)=\CM_{D_2}(D_2)=0$, so Theorem \ref{thm:cmNotEqualToCmSd} does not apply.  In the next section, we explore what $\CMD=0$ tells us about $D$.

\subsection{Consequence of $\CMD=0$}

In this section, we focus exclusively on properties of $\CMD$.  
Since  the results do not involve $\CM(D),$  we no longer require the assumption that $D$ is a finite set.  
However, for the sake of simplicity, we continue to assume that $D$ is a nonempty  set of  nonzero vectors. 

To explore the consequences of $\CMD=0$ we  first begin by examining the other extreme case: $\CMD=1$.   

\begin{proposition}\label{prop:CMD=1} Let $D \subseteq \R^n$ be a nonempty closed set of nonzero vectors.  
  Define $U=\{u \in \spann(D): \Vert u \Vert=1\}$ and $\widehat{D} = \{d/\|d\| : d \in D\}$.  Then the following are equivalent.
\begin{enumerate}[(i)]
\item $U = \widehat{D}$,
\item $\CM_D(D)=1$, 
\item $\V_D(D)=U$, 
\item $\A_D(D)=U$.
\end{enumerate}
\end{proposition}
\begin{proof} Clearly (i) implies (ii), (iii), and (iv).   Therefore we only need to show (ii) implies (i).  \\
Suppose that $\CM_{D}(D)=1$.  By construction, $\widehat{D} \subseteq U$.  
Let $u^* \in U$.  By definition of the cosine measure relative to $\spann(D)$,
    $$1 = \CM_D(D) = \min_{\substack{u \, \in \, \spann(D)\\\Vert u \Vert=1}} \sup_{d \, \in D} \frac{u^\top d}{\Vert d \Vert} =  \min_{\substack{u \, \in \, \spann(D)\\\Vert u \Vert=1}} \sup_{\hat{d} \in \widehat{D}} u^\top \hat{d} \leq \sup_{\hat{d} \in \widehat{D}} (u^*)^\top \hat{d} \leq \|u^*\| = 1,$$
where the last inequality is true by Cauchy–Schwarz.  Moreover, the last inequality is an equality if and only if $\hat{d}=u^*$.  Since equality holds across the above, this implies $u^* \in \widehat{D}$.  Thus $U \subseteq \widehat{D}$, which provides (i) and the proof is complete.\qed
\end{proof}

An important consequence of Proposition \ref{prop:CMD=1} is that if $|D|=2$, then either $\CMD=1$ or $\CMD<0$. 

\begin{corollary}\label{cor:Dhas2elements}
    Let $D \subseteq \R^n$ be a set of exactly $2$ nonzero vectors.  Then either $\CMD=1$ or $\CMD<0$. 
\end{corollary}
\begin{proof}Let $D=\{d_1, d_2\}$ and set $u = -\left(\frac{d_1}{\|d_1\|}+\frac{d_2}{\|d_2\|}\right)$.  If $u = \zero_n$, then Proposition \ref{prop:CMD=1} creates $\widehat{D}=\{\frac{d_1}{\|d_1\|}, -\frac{d_1}{\|d_1\|}\}=U$, so $\CMD=1$.  If $u \neq \zero_n$, then set $u^* = u/\|u\|$ and notice
    $$\CMD \leq \max\left\{ \frac{(u^*)^\top d_1}{\|d_1\|}, \frac{(u^*)^\top d_2}{\|d_2\|} \right\} = \frac{1}{\|u\|}\left(-1 - \frac{d_1^\top d_2}{\|d_1\|\|d_2\|}\right) < 0.$$ \qed
\end{proof}

We now begin our examination of the consequences of $\CMD=0$.  Our goal is to show that $\CMD=0$ if and only if $D$ contains a finite nontrivial subset that is a positive spanning.  The proof is reductionist in nature and uses the following lemma, which shows that if $\CMD=0$, then there exists a subset $V \subseteq D$ that is strictly smaller than $D$ and has $\CM_{V}(V)\geq 0$.  We assume that the set is closed  in the next  results.  Hence, the set considered can be either a finite set or a closed infinite set.    

\begin{lemma}\label{lem:reducingD}
    Let $D \subseteq \R^n$ be a closed set of nonzero vectors.  If $\CMD=0$, then there exists a finite subset $V \subseteq D$ such that $1 < |V| < |D|$ and $\CM_{V}(V) \geq 0$. 
\end{lemma}

\begin{proof}
Let $u^* \in \CVD$.  Define $H =\{v \in \spann(D) : (u^*)^\top v = 0\}$ and $V = D \cap H$.  Since $\CMD=0$ and $D$ is closed, by definition there exists at least one vector in $V$.

If $|D|=|V|$, then $V=D$, which implies $(u^*)^\top d = 0$ for all $d \in D$.  This yields $(u^*)^\top v = 0$ for all $v \in \spann(D)$, which implies $u^* = \zero_n$ contradicting $\|u^*\|=1$.  Thus, $|V|<|D|$.

Finally, for eventual contradiction, suppose that $\CM_{V}(V) < 0$.  This implies that there exists $v^* \in \spann(V)$ such that   $(v^*)^\top d < 0$ for all $d \in V$.  For $\epsilon >0$, define 
    $$u' = u^* + \epsilon v^* \in \spann(D).$$
Given any $d \in V$, we have 
    \begin{align*}
    (u')^\top d &= (u^* +  \epsilon v^*)^\top d
    = (u^*)^\top d + ( \epsilon v^*)^\top d < 0,
\end{align*}
as $(u^*)^\top d =0$ and $( \epsilon v^*)^\top d < 0$.  Given $d \in D \setminus V$, we have 
\begin{align*}
     (u')^\top d &= (u^*)^\top d + (\epsilon v^*)^\top d 
    \leq \max_{d \in D\setminus V}\left(\|d\|\frac{(u^*)^\top d}{\|d\|}\right)+ (\epsilon v^*)^\top d.
\end{align*}
Since $u^* \in \CVD$, we have $$\max_{d \in D \setminus V}\left(\|d\|\frac{(u^*)^\top d}{\|d\|}\right) < 0.$$ Thus, for $\epsilon$ sufficiently small $(u')^\top d < 0$ for all $d \in D$.  The vector $u'$ contradicts $\CMD = 0$, and therefore we must have $\CM_{V}(V) \geq 0$.

Finally, $|V|>1$, as  if $|V|=1$, then $\CM_{V}(V) = -1$.  If $V$ is not finite, then apply Theorem \ref{thm:equivalencepss} to reduce to a finite subset.
    \qed
\end{proof}

We now present the main result for this subsection.

\begin{theorem} \label{thm:containsProperPSS}
Let $D \subseteq \R^n$ be a nonempty closed set of nonzero vectors. Suppose that $\pspan(D) \neq \spann(D)$.  Then $\CMD=0$  if and only if $D$ contains a nonempty finite proper subset $V$ such that $\pspan(V) = \spann(V)$.
\end{theorem}
\begin{proof} First, note that  $\pspan(D) \neq \spann(D)$ implies $\CMD \leq 0$.

($\Leftarrow$)  Suppose that there exists a nonempty proper subset $V \subset D$ such that $\pspan(V)=\spann(V).$ Let $u^\ast \in \CVD$ and define $u^\ast_V = \ProjV u^\ast$ and  $u^\ast_{V^\perp} = u^\ast - u^\ast_V$.  Since $u^\ast_{V^\perp}$ is in the orthogonal subspace to $V$, we have 
    \begin{equation}\label{eq:CMinequalitysubspace}
   \sup_{d \, \in \, V} \frac{(u^\ast)^\top d}{\Vert d \Vert} =  \sup_{d\, \in \, V} \frac{(u_V^\ast)^\top d}{\Vert d \Vert} \leq  \sup_{d \, \in \, D} \frac{(u^\ast)^\top d}{\Vert d \Vert}= \CMD 
   \leq 0. \end{equation}
If $u^\ast_V \neq \zero_n,$ then Proposition \ref{prop:dotProdPos} would imply the existence of $d \in V$ with $(u_V^\ast)^\top d > 0$, therefore equation \eqref{eq:CMinequalitysubspace} implies that $u^\ast_V = \zero_n$.  Substituting $u^\ast_V = \zero_n$ in equation \eqref{eq:CMinequalitysubspace} shows $\CMD=0$.

($\Rightarrow$) Suppose that $\CMD=0$.  Without loss of generality, we assume $D$ is finite.  (Indeed, if $D$ is not finite, then apply Theorem \ref{thm:equivalencepss} to drop to a finite set.)  Let $m = |D|$.  Note that $m \geq 3$, as $|D|=1$ implies $\CMD=-1$ and $|D|=2$ implies $\CMD = 1$ or $\CMD<0$.

By Lemma \ref{lem:reducingD}, there exists a nonempty proper subset $D_{1} \subseteq D$ such that $1 <|D_{1}| < m$ and $\CM_{D_1}(D_1) \geq 0$.  If $\CM_{D_{1}}(D_{1}) > 0$, then $V = D_{1}$ is our desired set.  If $\CM_{D_{1}}(D_{1}) = 0$, then $|D_1| > 2$ and therefore 
we can repeat the procedure as necessary to generate a nonempty proper subset $D_{k} \subset D_{k-1}$ such that $2 \leq  |D_{k}| < |D_{k-1}| \leq m - k$ and $\CM_{{D_{k}}}(D_{k}) \geq 0$.  This process must terminate before $k = m - 1$ or a contradiction is created.  When the procedure is terminated we have $\CM_{{D_{k}}}(D_{k}) > 0$, so $V=D_k$ is our desired set. \qed
\end{proof}

The previous theorem can be adapted to the cosine measure of $D$ or reformulated to discuss dimensions of subspaces.

\begin{corollary}
Let $D \subseteq \R^n$ be a nonempty closed set of nonzero vectors.
\begin{enumerate}[(i)]
\item Suppose that $\pspan(D) \neq \R^n$.  Then $\CM(D)=0$  if and only if  $D$ contains a nonempty proper subset $V$ such that $\pspan(V) = \spann(V)$.
\item Suppose that $\pspan(D) \neq \spann(D)$.  Then $\CMD=0$  if and only if  $D$ contains a positive basis of a linear subspace $L \subset \spann(D)$ with $1 \leq \dim(L) <\dim(\spann(D)).$
\item Suppose that $\pspan(D) \neq \spann(D)$.  Then $\CMD=0$ if and only if  $D$ contains a minimal positive basis of a linear subspace $L \subset \spann(D)$ with $1 \leq \dim(L) <\dim(\spann(D)).$
\end{enumerate}
\end{corollary}

\begin{proof}
    Item (i) is immediate from Theorem \ref{thm:containsProperPSS}.  Item (ii) results from rephrasing Theorem \ref{thm:containsProperPSS} in terms of positive bases.  Item (iii) follows from  \cite[Theorem 1]{Romanowicz1987}, where it is shown that a positive basis of a linear subspace can be partitioned to minimal positive bases (see also \cite{Davis1954}). \qed
\end{proof}

\subsection{Bounding the norm of the gradient  in directional direct-search methods}

Directional methods such as the Generalized Pattern Search (GPS)~\cite{Torczon1997}
    and the Mesh Adaptive Direct Search (MADS)\cite{Audet2006} algorithms
    are foundational methods in derivative-free optimization \cite{audet2017derivative,conn2009introduction}.  
Their convergence rely on applying poll steps, evaluating $f(x^k + d)$, $d \in D$ where $D$ is a positive spanning set for $\R^n$, to seek improvement in the objective function.   An important result is that in the case of a failed poll step (i.e., $f(x^k) \leq f(x^k + d)$ for all $d \in D$), the cosine measure of $D$ combined with the radius of $D$ provides an error bound on the gradient of $f$ (assuming  $\nabla f$ is Lipschitz continuous).  We provide a formal statement below from \cite{conn2009introduction}, but recommend seeing \cite{dolan2003,kolda2003} for alternate presentations.  

Given the nature of directional direct-search algorithm, in this subsection we continue to assume that $\zero_n \notin D$.

\begin{theorem}\cite[Theorem 2.8]{conn2009introduction} \label{thm:gradientBound} 
Let $D \subseteq \R^n$ be a nonempty finite subset of nonzero vectors with radius $\Delta_D$.  Let $f:\dom f \subseteq \R^n \to \R$ and $x^0 \in \dom f$.  
Suppose that $D$ positively spans $\R^n$ and $f(x^0) \leq f(x^0+d)$ for all $d \in D.$  If $\nabla f$ is Lipschitz continuous with constant $L_{\nabla f} \geq 0$ in an open set containing the ball $B_n(x^0;\Delta_D)$, then 
\begin{align*}
    \Vert \nabla f(x^0) \Vert &\leq \frac{1}{2} L_{\nabla f} \CM(D)^{-1} \Delta_D.
\end{align*}
\end{theorem}

We next present two extensions of this result.  Theorem \ref{thm:EBPSSD}(i) directly extends Theorem \ref{thm:gradientBound} to allow for subspaces and infinite sets.  (Setting $L=\R^n$ and $D$ finite in Theorem \ref{thm:EBPSSD}(i) reproduces Theorem \ref{thm:gradientBound}.)  Theorem \ref{thm:EBPSSD}(ii) presents a new result, demonstrating a stronger error bound in the situation where the set $D$ contains symmetry. 

\begin{theorem} \label{thm:EBPSSD} 
Let $D \subseteq \R^n$ be a nonempty (possibly infinite) set of nonzero vectors.  
Suppose that the radius of $D$ is finite $(0<\Delta_D<\infty)$.  Let $f:\dom f \subseteq \R^n \to \R$ and $x^0 \in \dom f$.  
In addition, suppose that $\pspan(D)=\spann(D)$ and $f(x^0) \leq f(x^0+d)$ for all $d \in D.$
\begin{enumerate}[(i)]
\item  If $\nabla f$ is Lipschitz continuous with constant $L_{\nabla f} \geq 0$ in an open set containing the ball $\overline{B}_n(x^0;\Delta_D)$, then 
\begin{align} \label{eq:secondBound}
    \Vert \ProjD \nabla f(x^0) \Vert  &\leq \frac{1}{2} L_{\nabla f} \CMD^{-1} \Delta_D.
\end{align}
\item  Suppose in addition that for each $d\in D$ there exists $\alpha_d > 0$ such that $-\alpha_d d \in D$.  If $\nabla^2 f$ is Lipschitz continuous with constant $L_{\nabla^2 f} \geq 0$ in an open set containing the ball $\overline{B}_n(x^0;\Delta_D)$, then 
\begin{align} \label{eq:thirdBound}
    \Vert \ProjD \nabla f(x^0) \Vert  &\leq \frac{1}{3} \alpha_{\max}L_{\nabla^2 f}\CMD^{-1}  \Delta_D^2, 
\end{align}
where $\alpha_{\max} =\sup_{d \in D}\{\alpha_d\}$.
\end{enumerate}
\end{theorem}
\begin{proof} (i)  Let $v=-\ProjD \nabla f(x^0)$.  If $v \neq \zero_n$, then by definition, we have
\begin{align*}
    \CMD &\leq \sup_{d \in D} \frac{v^\top d}{\Vert v \Vert \Vert d \Vert}
\end{align*}
and therefore there exists a $d \in D$ such that 
\begin{align*}
    \CMD\Vert v \Vert \Vert d \Vert &\leq v^\top d.
\end{align*}
If $v = \zero_n$, the above holds trivially.  Therefore, there exists a vector $d \in D$ such that 
\begin{align} \label{eq:firstIneq}
 \CMD \Vert \ProjD \nabla f(x^0) \Vert \Vert d \Vert  &\leq -\left (\ProjD \nabla f(x^0) \right )^\top d.
\end{align}
Applying Taylor's Theorem and the assumption that $f(x^0)-f(x^0+d) \leq 0$ for all $d \in D$, we have
\begin{align}\label{eq:secondIneq}
    -\left (\ProjD \nabla f (x^0) \right )^\top d = -\nabla f(x^0)^\top d&=f(x^0)-f(x^0+d)+R_1(x^0;d) \leq  R_1(x^0;d)
\end{align}
where $R_1$ is the first-order remainder term. Combining equations \eqref{eq:firstIneq} and \eqref{eq:secondIneq}, then noting that $|R_1(x^0;d)| \leq \frac{1}{2} L_{\nabla f}\|d\|^2$, yields
\begin{align*}
     \CMD \Vert \ProjD \nabla f(x^0) \Vert \Vert d \Vert  &\leq   \vert R_1(x^0;d) \vert \leq \frac{1}{2} L_{\nabla f} \Vert d \Vert^2.
\end{align*}
Since $\spann(D)=\pspan(D)$, we know that $\CMD>0.$ Equation \eqref{eq:secondBound} follows using the fact $\Vert d \Vert \leq \Delta_D$ for all $d \in D.$

(ii) Suppose that for each $d\in D$ there exists $\alpha_d > 0$ such that $-\alpha_d d \in D$ and let $\alpha_{\max} = \sup_{d \in D}\{\alpha_d\}$. If $\alpha_{\max} = \infty$, then the result holds trivially, so we assume $\alpha_{\max} < \infty$. 

Taylor's Theorem and the assumption that $f(x^0)-f(x^0+d) \leq 0$ for all $d \in D$, yields 
\begin{align} 
    -\nabla f(x^0)^\top d&=f(x^0)-f(x^0+d)+\frac{1}{2} d^\top \nabla^2 f(x^0) d+R_2(x^0;d),\nonumber\\
    -\nabla f(x^0)^\top d&\leq \frac{1}{2} d^\top \nabla^2 f(x) d+R_2(x^0;d)\label{eq:Taylor2d}
\end{align}
and similarly,
\begin{align}
 -\nabla f(x^0)^\top (-\alpha_d d)&\leq \frac{\alpha_d^2}{2} d^\top \nabla^2 f(x^0) d+R_2(x^0;-\alpha_d d),   \label{eq:Taylor2Minusd}
\end{align}
where $R_2$ is the second-order remainder term.  
Multiplying \eqref{eq:Taylor2d} by $\alpha^2$ and subtracting  \eqref{eq:Taylor2Minusd}, we get
\begin{align} \label{eq:boundOrder2}
    -(\alpha_d^2+\alpha_d)\nabla f(x^0)^\top d= -(\alpha_d^2+\alpha_d) \left (\ProjD \nabla f(x^0) \right )^\top d 
    \leq \alpha_d^2 R_2(x^0;d)-R_2(x^0;-\alpha_d d)
\end{align}
Therefore, combining equations \eqref{eq:firstIneq} and \eqref{eq:boundOrder2},  then noting that $|R_2(x^0;d)| \leq \frac{1}{6} L_{\nabla^2 f}\|d\|^3$, we find
\begin{align*}
   (\alpha_d^2+\alpha_d) \CMD \Vert \ProjD \nabla f(x^0) \Vert \Vert d \Vert  &\leq  \alpha_d^2 R_2(x^0;d)-R_2(x^0;-\alpha_d d) \\
    &\leq \alpha_d^2 \frac{L_{\nabla^2 f}}{6} \Vert d \Vert^3 + \frac{L_{\nabla^2 f}}{6} \Vert -\alpha_d d \Vert^3 \\
    &\leq (\alpha_d^2 +\alpha_d^3) \frac{L_{\nabla^2 f}}{3} \Vert d \Vert^3
\end{align*}
Using the bounds $\Vert d \Vert \leq \Delta_D$ and $\alpha_d \leq \alpha_{\max}$, we obtain \eqref{eq:thirdBound}. 
\qed
\end{proof}

Theorem \ref{thm:EBPSSD} could be used in directional direct-search methods in several different manners.  
One obvious example would be to use it to generate a stopping condition, particularly in the case where $\spann(D)=\pspan(D)$.  
In the case of reduced subspace methods where $\spann(D) \neq \pspan(D)$, Theorem \ref{thm:EBPSSD} could be used to create flags indicating when it is time switch to a different subspace.

While convergence of the directional direct-search methods requires the use of positive spanning sets, it is easy to conceive of an implementation that does not enforce positive spanning sets at every iteration.  The following corollary demonstrates how Theorem \ref{thm:EBPSSD} might be used to help determine next steps in the case of a failed poll step where $D$ is not a positive spanning set.

\begin{corollary} \label{cor:EBNotPSSSD} Let $D \subseteq \R^n$ be a nonempty (possibly infinite)  set of nonzero vectors.  Suppose that the radius of $D$ is finite $(0<\Delta_D<\infty)$.  Let $f:\dom f \subseteq \R^n \to \R$ and $x^0 \in \dom f$.  
In addition, suppose that $\pspan(D) \neq \spann(D)$ and $f(x^0) \leq f(x^0+d)$ for all $d \in D.$

Let $B=\{d_1, d_2 \ldots, d_m\} \subseteq D$ be  a basis of $\spann(D)$.  Define  $$w=-\Delta_D \frac{\sum_{j=1}^m d_j}{\|\sum_{j=1}^m d_j\|}, \quad D'=D \cup \{w\}, \quad \mbox{and} \quad D''=D \cup -D.$$

\begin{enumerate}[(i)]
\item If $\nabla  f$ is Lipschitz continuous  with constant $L_{\nabla f} \geq 0$ in an open set containing the ball $\overline{B}_n(x^0;\Delta_D)$, then at least one of the following holds:
\begin{align}f(x^0+w) < f(x^0) \quad\mbox{\em or}\quad
    \Vert \ProjD \nabla f(x^0) \Vert  &\leq \frac{1}{2} L_{\nabla f} \CM_{D}(D')^{-1} \Delta_{D}.
\end{align}

\item If $\nabla^2 f$ is Lipschitz continuous  with constant $L_{\nabla^2 f} \geq 0$ in an open set containing the ball $\overline{B}_n(x^0;\Delta_D)$, then then at least one of the following holds
\begin{align}f(x^0-d) < f(x^0) \mbox{ for some } d \in D \quad\mbox{\em or}\quad
    \Vert \ProjD \nabla f(x^0) \Vert  &\leq \frac{1}{3} L_{\nabla^2 f}\CM_{D}(D'')^{-1} \Delta_D^2. 
\end{align}
\end{enumerate}
 \end{corollary}
\begin{proof}(i) Suppose that $f(x^0)\leq f(x^0+w)$.  By Proposition \ref{prop:extendPSSSD}, $D'$ is a positive spanning set.  The result follows from Theorem \ref{thm:EBPSSD}(i), noting that the scaling of $w$ makes $\Delta_{D'}=\Delta_D$.  
\\
(ii) Suppose that $f(x^0) \leq f(x^0 - d)$ for all $d \in D$. By Proposition \ref{prop:extendPSSSD}, $D''$ is a positive spanning set.  The result follows from Theorem \ref{thm:EBPSSD}, noting that $\Delta_{D''}=\Delta_D$ and $\alpha_d=1$ for all $d \in D$. \qed
\end{proof}

\section{Computing the cosine measure relative to a subspace} \label{sec:algo}

Section \ref{sec:CMrelativetoL} established the value of the cosine measure relative to a subspace. In this section, we investigate how to compute $\CMD$ for a nonempty finite set of nonzero vectors. In Example \ref{ex:valueCMrelative}, we were able to compute the cosine measure relative to a subspace since the subspace considered  were either 1-dimensional or 2-dimensional.  In this section, we provide  a general deterministic algorithm to compute the cosine measure relative to a subspace. Note that Algorithm \ref{alg:cmpspanning} assumes $D$ is finite and $\zero_n \notin D$.  Reasons for both are clear.  If we allowed $\zero_n \in D$, then the first step of the algorithm would simply become remove $\zero_n$ from $D$. If $D$ is infinite, then we must confront the challenge of how to express the set in a manner suitable for algorithmic use. Developing an algorithm for general infinite sets is  an obvious future research direction to investigate.    

Algorithm \ref{alg:cmpspanning} is a modification of the algorithm in \cite{hare2020deterministic} to allow for subspaces and for the case where $D$ is not a positive spanning set.  To allow for the scenario where $D$ is not a positive spanning set, Algorithm \ref{alg:cmpspanning} begins by checking if $\pspan(D)=\spann(D)$. Recall that this can be done by solving a linear program using Theorem \ref{thm:equivalencepss}(vii). In this section, the notation $\G(B)$ represents the \emph{Gram matrix} of $B$. That is $\G(B)=B^\top B.$ 

\begin{center}
\begin{algorithm}
\DontPrintSemicolon
\caption{The cosine measure of a finite set $D$ relative to  $\spann(D)$ \label{alg:cmpspanning}}
Given a nonempty finite set $D$ of  $q$ nonzero vectors  in $\R^n,$

\textbf{0. Normalize:} set $D \leftarrow \{ d/\|d\| : d \in D\}$.

\textbf{1. Determine if $\pspan(D)=\spann(D).$}

   \quad  If $\pspan(D) =\spann(D),$ go to 2.\\
    \quad If $\pspan(D) \neq \spann(D),$ go to 3.

\textbf{2. Let $m=\dim(\spann(D))\geq 1.$ For all bases $B$ of $\spann(D)$ contained in $D,$  compute\;}

$\begin{aligned}
    \quad (2.1)&\quad\gamma_B= \frac{1}{\sqrt{\one_m \top \G(B)^{-1}\one_m}} &&\text{(positive value of the $m$ equal dot products)},\\
    \quad (2.2)&\quad u_B=\gamma_{B} \left (B^{\top} \right )^\dagger  \one_m &&\text{(the unit vector associated to $\gamma_{B}$)},\\
    \quad (2.3)&\quad p_B=\begin{bmatrix}[p_{B}]_1&\cdots&[p_{B}]_q \end{bmatrix}=(u_B)^\top D &&\text{(the dot product vector)}, \\
    \quad (2.4) & \quad\p =\max_{1\leq j \leq q }[p_B]_j &&\mbox{(the maximum value in $p_{B}$)}.\\
\end{aligned}$

\textbf{3. Return \:} \\
$\begin{aligned}
    \quad (3.1)&\quad \CMD=\left\{\begin{array}{lll}
         \displaystyle\min_{B \, \subseteq \, D} \p, & \mbox{if $\pspan(D) =\spann(D)$},\\
         \displaystyle \min_{\substack{u \, \in \, \spann(D)\\\Vert u \Vert \leq 1}}\max_{d \in D} u^\top d, & \mbox{if $\pspan(D) \neq \spann(D)$}.\\
   \end{array}\right.\\
    \quad (3.2)&\quad\CVD=
        \left\{\begin{array}{lll}
        \displaystyle\{u_{B}:\p=\CMD\}, &  \mbox{if $\pspan(D)=\spann(D)$}, \\
        \displaystyle \{u \in \spann(D):u^\top D\leq \zero_q, \, \Vert u \Vert=1\} & \mbox{if $\pspan(D) \neq \spann(D).$}\\
   \end{array}\right.\\   
\end{aligned}$
\end{algorithm}
\end{center}

Algorithm \ref{alg:cmpspanning} seperates the computation of $\CMD$ into two scenarios. The first scenario, $\pspan(D) = \spann(D)$, is solved through a finite enumeration process similar to \cite{hare2018computing}. 
The second scenario,  $\pspan(D) \neq \spann(D)$, is solved by applying the convex relaxation to the nonconvex cosine measure problem.  Notice that the convex relaxation is a second-order conic program, so it can be solved by a variety of well-studied algorithms.  Our implementation uses the function \emph{coneprog} in MATLAB.

The next lemma shows that if $\pspan(D) \neq \spann(D)$, then the convex relaxation of the cosine measure problem can be considered.

\begin{lemma} \label{lem:convexrelaxexact}
Let $D$ be a nonempty finite set of nonzero vectors in $\R^n$ with $\pspan(D) \neq \spann(D)$.  Then the convex relaxation of the cosine measure problem is exact; i.e., 
    \begin{equation}\label{eq:convexrelaxexact}
        \CMD = \min_{\substack{u \, \in \, \spann(D)\\\Vert u \Vert \leq 1}}\max_{d \in D} u^\top d.
    \end{equation}
Moreover, if $\CMD \neq 0$, then
    $$\CVD = \argmin{\substack{u \, \in \, \spann(D)\\\Vert u \Vert \leq 1}}\max_{d \in D} u^\top d.$$
\end{lemma}

\begin{proof}
Denote 
    $$p^* = \min_{\substack{u \, \in \, \spann(D)\\\Vert u \Vert \leq 1}}\max_{d \in D} u^\top d \quad \mbox{and} \quad X^* = \argmin{\substack{u \, \in \, \spann(D)\\\Vert u \Vert \leq 1}}\max_{d \in D} u^\top d.$$
It suffices to show that $X^*$ contains at least one vector $x^*$ with $\|x^*\|=1$.

We claim that $X^*\setminus\{\zero\} \neq \emptyset$. As $D$ is finite, $X^*$ is nonempty.  If $X^*=\{\zero\}$, then $p^*=0$ and for any $u \in \spann(D)$ with  $0 <\Vert u \Vert \leq 1,$  we have $\max_{d \in D} u^\top d > 0$.  This would imply that $\CMD > 0$, which contradicts Proposition \ref{prop:pssSDiffcmGreaterEqualZero}.  Hence $X^*\setminus\{\zero\} \neq \emptyset$.

Select any $x^* \in X^*\setminus\{\zero\}$.  Since $p^* \leq \CMD \leq 0$, we have $\max_{d \in D} (x^*)^\top d \leq 0$.  Using $\|x^*\|\leq 1$, this implies $\max_{d \in D} \frac{x^*}{\|x^*\|}^\top d \leq \max_{d \in D} (x^*)^\top d = p^*$, so $\frac{x^*}{\|x^*\|} \in X^*$.  As such,  $X^*$ contains at least one point $x^*$ with $\|x^*\|=1$, and $\CMD=p^*$ as desired.  The final statement follows trivially; if $\CMD \neq 0$, then $ \zero \notin X^*$. \qed
\end{proof}

\begin{remark}
    When $\pspan(D) \neq \spann(D)$, Lemma \ref{lem:convexrelaxexact} can always be used to find $\CMD$.  However, when $\CMD = 0$, some second-order conic program solvers may return $0$ as the minimizer for problem \eqref{eq:convexrelaxexact}.  In this case, an orthonormal basis of the null space  of the sub-positive spanning set contained in $D,$ say $B_{ns},$ may  be found (for example using the MATLAB function \emph{null}).  A cosine vector is obtained by verifying  the inequality  $u^\top D\leq \zero_q$  for any vector  $u \in B_{ns} \cup (-B_{ns})$.  
\end{remark}

To prove that the algorithm returns the correct cosine measure and the cosine vector set in the case where $\pspan(D)=\spann(D)$,  we begin by introducing  several results that can be viewed as  a generalization of the results in \cite{hare2020deterministic,naevdal2018}. 
The following lemma is an adaptation   of Lemma 1 in \cite{naevdal2018}.The proof is essentially identical to the proof of \cite[Lemma 1]{naevdal2018} by replacing $\R^n$ with $S_B.$
\begin{lemma} \label{lem:gammab}
    Let $B=\begin{bmatrix} d_1&\cdots&d_m \end{bmatrix}$ be a linearly independent set of vectors  in $\R^n$ written in matrix form where $1 \leq m \leq n$ and where each $d_j$ is a unit vector.
    Then there exist a unit vector $u_B \in S_B=\spann(B)$  such that   $$(u_B)^\top d_1=\dots=(u_B)^\top d_m=\gamma_B,$$
    and $$((-1)u_B)^\top d_1=\dots=((-1)u_B)^\top d_m=(-1)\gamma_B.$$
    where 
    \begin{align} 
    \gamma_B&=\frac{1}{\sqrt{\one_m^\top  \G(B)^{-1}\one_m}}>0. \label{eq:gammaBplus}
    \end{align}
    Moreover,  
    \begin{align} \label{eq:vectorUbplus}
    u_B&= \gamma_{B} (B^\top)^\dagger  \one_m.
    \end{align}
\end{lemma}

\begin{lemma} \label{lem:alphaEqualgammab}
Let $B=\{d_1, \cdots, d_m\}$ be a linearly independent set of vectors in $\R^n$ where each $d_j$ is a unit vector. Suppose that $u$ is a unit vector in $S_B=\spann(B)$ such that $u^\top d_1=\dots=u^\top d_m=\alpha>0.$ Then $\alpha=\gamma_B,$ where $\gamma_B$ is defined as in \eqref{eq:gammaBplus}.
\end{lemma}
The previous lemma can be proved using a similar process than the proof of Lemma 13 in \cite{hare2020deterministic}. Next, we provide lemma that will be useful to prove the key theorem in the case $\pspan(D)=\spann(D).$

\begin{lemma}\label{lem:epsover2}
Let $\epsilon \neq 0$ and let $u$ and $v$ be unit vectors in $\R^n.$  Then
    \begin{enumerate}[(i)]
        \item $\Vert u+\epsilon v \Vert=1$ if and only if $\epsilon=-2u^\top v$, and
        \item $\Vert u+\epsilon v \Vert<1$ implies $\Vert u-\epsilon v \Vert>1$.
        \item Assume   $\Vert u \pm \epsilon v \Vert \neq 0.$ Then 
        \begin{align*}
          \frac{u \pm \epsilon v}{\Vert u \pm \epsilon v \Vert}= \pm u \iff v =\pm u.
        \end{align*}
        
    \end{enumerate}
\end{lemma}
\begin{proof}
Since
   $$
    \Vert u +\epsilon v \Vert^2=1 + (2\epsilon u^\top v + \epsilon^2) ~\mbox{and}~\Vert u -\epsilon v \Vert^2=1 - (2\epsilon u^\top v - \epsilon^2),
   $$
it follows that $\Vert u +\epsilon v \Vert=1$ if and only if $2\epsilon u^\top v + \epsilon^2=0$.  Since $\epsilon \neq 0$, the first result follows.

The second result follows by contradiction.
Suppose that $\|u+\epsilon v\|^2 < 1$ and
$\|u-\epsilon v\|^2 \leq 1$.
By expanding and adding these inequalities, we get the contradiction $2\epsilon^2 < 0$.

To show the third result,
 suppose that $\Vert u \pm \epsilon v \Vert \neq 0$ and
 observe that 
\begin{eqnarray*}
    \frac{u \pm \epsilon v}{\Vert u \pm \epsilon v \Vert} = \pm u
    & \Leftrightarrow & 
    v = \lambda u \qquad \mbox{ where } \lambda = \frac1{\pm\epsilon}\left( \pm 1 - \|u \pm \epsilon v\| \right).
\end{eqnarray*}
But since both $u$ and $v$ are unit vectors, it follows that the last equality is equivalent to $v = \pm u$.

\qed
\end{proof}

The next theorem contains the main result  governing Algorithm \ref{alg:cmpspanning} in the case $\pspan(D)=\spann(D)$. It assumes that the set $D$ is finite. However,  note that the following two results (Theorem \ref{thm:spanSD} and Corollary \ref{cor:containbasis}) hold for closed infinite sets. But since the algorithm  assumes $D$ is finite, we keep the same assumption and assume $D$ is finite.
 \begin{theorem}\label{thm:spanSD}
Let $D$ be a nonempty  finite set of nonzero vectors in  $\R^n$   such that $\pspan(D)=\spann(D).$  Let $\us \in \CVD$.
Then $$\spann(\A_D(D,\us))=\spann(D).$$
\end{theorem}
\begin{proof} Without loss of generality, assume that all vectors $d$ in $D$ are unit vectors. Clearly we have $\spann(\A_D(D,u_B)) \subseteq \spann(D)$ since $\A_D(D,u_B) \subseteq D.$ It remains to show that $\spann(D) \subseteq \spann(\A_D(D,u_B)).$ 
Suppose that $\spann(D) \nsubseteq \spann(\A_D(D,\us))$. This means that the rank of $\A_D(D,\us)$ is strictly less than $\dim(\spann(D)).$ Since $D$ is closed, we obtain  $$1 \leq \dim(\spann(\A_D(D,u_B))<\dim(\spann(D)).$$  From the previous inequalities and $\spann(\A_D(D,u_B) \subseteq \spann(D),$ it follows that there exists   a unit  vector $v \in \spann(D)$ which is in the kernel of $\A_D(D,\us)$. This means that $d^\top v=0$ for all $d$ in $\A_D(D,\us).$

Notice that, if $d \in D\setminus \A_D(D,\us),$ then $$d^\top \us <\CMD.$$
 Consider the vector $\us+\epsilon v \in \spann(D).$  Since $\pspan(D)=\spann(D),$ we must have $\CMD>0,$ and it follows that  $\us \neq \pm v$ as $\us^\top d \neq 0$ for  all $d \in \A_D(D,\us).$   Hence, using Lemma \ref{lem:epsover2}(iii),  for sufficiently small $\epsilon>0$ and   not equal to $\vert -2\us^\top v\vert,$  we have  
    $$\frac{d^\top (\us\pm \epsilon v)}{\Vert \us \pm \epsilon v \Vert}<\CMD$$
for all $d \in D\setminus \A_D(D,\us).$  Moreover, since $d^\top v=0$  for all $d \in \A_D(D,\us),$ it follows that 
    $$\frac{d^\top(\us\pm \epsilon v)}{\Vert\us \pm \epsilon v \Vert}=\frac{d^\top \us}{\Vert\us \pm \epsilon v \Vert} \pm 0=\frac{\CMD}{\Vert \us \pm \epsilon v \Vert}$$
for all $d \in \A_D(D,\us)$ and where $\CMD>0.$
By Lemma \ref{lem:epsover2}(i), $\epsilon \neq -2\us^\top v$ implies that $\Vert \us + \epsilon v \Vert \neq 1.$ By Lemma \ref{lem:epsover2}(ii), if $\Vert \us +\epsilon v \Vert<1$, then $\Vert \us -\epsilon v \Vert>1.$ Select $w$  in $\{\us +\epsilon v, \us -\epsilon v \}$ such that $\Vert w \Vert >1.$ Then 
    $$\frac{d^\top w}{\Vert w \Vert}< \CMD$$ 
for all $d \in D.$ This contradicts the definition of cosine measure.

Therefore, $\spann(D)\subseteq \spann(\A_D(D,\us))\subseteq \spann(D)$, and the result follows. \qed
\end{proof}

Theorem \ref{thm:spanSD} can be viewed as an extension of Proposition 17 in \cite{hare2020deterministic} for the following reason: the linear subspace $\spann(D)$ is considered rather than the whole space $\R^n.$ 

The following corollary follows from Theorem \ref{thm:spanSD} and  the fact that a spanning set of a vector space  contains a basis of the vector space \cite[Theorem 2.11]{Brown1988}. 

\begin{corollary}\label{cor:containbasis}
Let $D$ be a nonempty finite set  of nonzero vectors  in $\R^n$  such that $\CMD>0.$ Let $\us \in \V(D).$ Then
$\A_D(D,\us)$ contains a basis of $\spann(D).$
\end{corollary}

We are now ready to show that Algorithm \ref{alg:cmpspanning} returns the desired values.

\begin{theorem}
Let $D=\{ d_1, \dots,d_q \}$ be   a  set  of $q \geq 1$ nonzero vectors in $\R^n.$ Then Algorithm \ref{alg:cmpspanning} returns $\CMD$ and $\CVD$.
\end{theorem}
\begin{proof}
Without loss of generality, assume that all vectors $d_j$ are unit vectors. Since $q\geq 1$, $\CVD \neq \emptyset$.   Consider the following two cases.

{\bf Case (i) $\pspan(D) \neq \spann(D).$} Then  Lemma \ref{lem:convexrelaxexact}  shows that Algorithm \ref{alg:cmpspanning} returns the value of the cosine measure and the cosine vector set.

{\bf Case (ii) $\pspan(D)=\spann(D).$} By Proposition \ref{prop:pssSDiffcmGreaterEqualZero}, this implies  $\CMD>0.$  Let $\us \in \CVD$. By Corollary \ref{cor:containbasis},  $\A_D(D,\us)$ contains a basis of $\spann(D)$.  Without loss of generality, let this basis (written in matrix form) be $B_{*}=\begin{bmatrix} d_1&\cdots&d_m \end{bmatrix}$ where  $m=\dim(\spann(D))\geq1.$ Hence,

\begin{align*}
    \CMD=d_1^\top \us=\cdots=d_m^\top \us>0.
\end{align*}
By Lemma  \ref{lem:alphaEqualgammab},
\begin{align*}
    \CMD=\gamma_{B_*}=\frac{1}{\sqrt{\one_m^\top \G(B_*)^{-1} \one_m}}
\end{align*}
Note that $\mathring{p}_{B_\ast}=\max_{1 \leq j \leq q} d_j^\top \us=\gamma_{B_{\ast}}$ since  $\gamma_{B_{\ast}}=\CMD.$ Therefore, we have 
\begin{align*}
    \CMD=\min_{B \subseteq D} \p=\mathring{p}_{B_{\ast}}.
\end{align*}
Taking all the  vectors $u_B$  associated to $\p$ such that $\CMD=\p$ in Step (3.2) returns the complete set $\CVD.$
 \qed
\end{proof}

\section{Conclusion} \label{sec:conclusion}
This paper introduces the definitions of the cosine measure relative to a subspace, and the related cosine vector set relative to a subspace.  
These definitions not only generalize the cosine measure and cosine vector set to allow for working in subspaces, but also generalize these ideas to work for infinite sets. 

When working with sets that are not positive spanning sets of $\R^n$, the novel definition of cosine measure  relative to a subspace may be  valuable as it provides more information about the set considered.
Proposition~\ref{prop:extendPSSSD} shows that any nonempty set of vectors can be extended to positive spanning set of its span by adding at most one vector to the set. 
Section~\ref{sec:CMrelativetoL} provides several properties of the cosine measure relative to a subspace. 
Theorem~\ref{thm:cmNotEqualToCmSd} shows that if the cosine measure of a set $D$ differs from the cosine measure relative to the span of $D$, then $D$ is  a positive spanning set of its span  and  $\spann(D)$ must be a  proper subspace of $\R^n.$  
Theorem~\ref{thm:containsProperPSS} proves that the cosine measure relative to $\spann(D)$ is equal to zero if and only if the set $D$ contains a positive spanning set of the span of a proper nonempty subset of $D$. 
Theorem \ref{thm:EBPSSD} uses the notion of cosine measure relative to a subspace to define two error bounds on the projected gradient of a smooth function. 
In the case where the set is not a positive spanning set of its span, Corollary~\ref{cor:EBNotPSSSD} introduces results that could be valuable when the poll step of a derivative-free algorithm fails.  
Lastly, a deterministic algorithm is proposed  to compute the cosine measure relative to its span. The algorithm is designed to accept  a non-positive spanning set as an input.  Combined, these results demonstrate that the cosine measure relative to a subspace is a valuable and practical tool to quantify the positive spanning properties of a set relative to a subspace. On a final note, an implementation of Algorithm \ref{alg:cmpspanning} in MATLAB is available upon request.


\vspace{1cm}
\small{\noindent \textbf{Acknowledgement} The authors  would like to thank the two referees for carefully reading our manuscript.  Their comments   have significantly improved our manuscript. The example showing that an infinite set can have an empty active set in Section 2 has been provided by a one of the referees.} 

\vspace{0.3cm}

\small{\noindent \textbf{Data availability}
No datasets were generated or analysed during the current study.}

\section*{Declarations}

\small{\textbf{Conflicts of interest} The authors have no competing interests to declare that are relevant to the content of this article.}

\normalsize

\bibliographystyle{siam}
\bibliography{bibliography}
\end{document}